\documentclass{elsart}
\usepackage{amsmath}
\usepackage{latexsym}
\usepackage{amssymb}
\newcommand{\n}{\infty}

\def\RR{\mathbb{R}}
\def\NN{\mathbb{N}}
\def\ZZ{\mathbb{Z}}
\def\CC{\mathbb{C}}

\begin{document}

\begin{frontmatter}

\title{ About composition of Toeplitz operators in Segal-Bargmann Space. }

\author[r]{Ramirez Romina,}
\author[m]{Marcela Sanmartino}

\address[r]{romina@mate.unlp.edu.ar}

\address[m]{mmasanmartino@gmail.com}
\address{Departamento de
Matem\'atica, Facultad de Ciencias Exactas,  Universidad Nacional
de La Plata, La Plata, Argentina}
\date{}

\newpage
\begin{abstract} Some positive results about the composition of Toeplitz operators
on the Segal-Bargmann space are presented. A Wick symbol
satisfying that it is not possible to construct its associated
Toeplitz operator is given.
\end{abstract}

%\noindent{\bf 2000 AMS Subject Classification}: Primary 47B35;
%Secondary 47G30
\end{frontmatter}
\newpage
\section{\textbf{Introduction}}
\label{intro}

 Toeplitz operators have been object of study in different
disciplines. In physics, these operators (also called anti Wick
operators) were introduced by Berezin as a quantization rule in
quantum mechanics (see \cite{B1}, \cite{B2}, \cite{S-B}). In
partial differential equations , Toeplitz operators and their
adjoints, play an important role in extending known results in the
space of entire functions to the context of Segal Bargmann
spaces(\cite{Ci}, \cite{Scha}), \cite{Ja1}, \cite{Ja} among
others). Also, they have been extensively studied as an important
mathematical tool in signal analysis (\cite{C-G}, \cite{C-R},
\cite{G}).

One of the problems still open is to know how to define the class
of the symbols where the composition of Toeplitz operators is
closed. In section 4, we give a positive result for radial
symbols.

Many authors have studied this problem. In  \cite{A-M},
\cite{C-G},  \cite{C-R} and \cite{L} , the composition of Toeplitz
operators  was shown to be written as a Toeplitz operator plus a
controllable remainder term, and in \cite{C} some exact results we
can be found.

\section{An overview of Toeplitz operators.}

Let $L^{2}(\mathbb{C}^{n}, d\mu)$ be the Hilbert space of square
integrable functions on $\mathbb{C}^{n}$, with the Gaussian
measure $d\mu(z)=\pi^{-n}e^{-|z|^{2}}dv(z)$ where $dv(z)$ is the
Lebesgue measure on $\mathbb{C}^{n}$.

In this work we deal with the \emph{Segal-Bargmann} space
$F^{2}(\CC^{n})$, the subspace of $L^{2}(\CC^{n}, d\mu)$ of all
the analytic functions in $\CC^{n}$. This space is a particular
case of a Fock space.

For more details and the motivation of this section, see the
Appendix.

Let $P:L^{2}(\CC^{n}, d\mu)\rightarrow {F}^{2}(\CC^{n})$ be the
orthogonal projection operator on the direction of the coherent
state $K_{w}(z)=e^{z.\overline{w}}$
\begin{equation}\label{projector}
Pg(z) = \langle g,K_{z} \rangle = \int_{\CC} g(w) \,
e^{z.\overline{w}} \, d\mu(w)\
\end{equation}
 where $\langle .  \ , . \rangle$ is the usual inner product in $L^{2}(\CC^{n},
d\mu)$.

Now, for a complex function $\varphi(z,\overline{z})$, we define
the operator  $T_{\varphi}$ on $F^{2}(\CC^{n}) $given by
 \begin{eqnarray}\label{toepl}
T_{\varphi}f(z)&:=& \int \varphi(w,
\overline{w}) \, e^{z\overline{w}} \, f(w) \, d\mu(w) \nonumber \\
    &=&P[\varphi f](z)\ .
\end{eqnarray}

Then, the domain of this operator is
$$\textit{Dom}(T_{\varphi})=\{f \ \mbox{such that} \  \varphi f \in
L^{2}(\CC^{n}, d\mu)\},$$ and $T_{\varphi}:
\textit{Dom}(T_{\varphi}) \rightarrow F^{2}(\CC^{n})$

\begin{defn}
\ \\
If for all $z \in \mathbb{C}^n$, $K_{z} \in
\textit{Dom}(T_{\varphi})$ we say that  $T_{\varphi}$ is a
\emph{Toeplitz operator} (\emph{or anti Wick operator}) with anti
Wick symbol $\sigma^{AW}(T_{\varphi})=\varphi$.
\end{defn}

\begin{rem}\label{domain}
\ \\
If $|\varphi(z)| \leq C e^{\delta r^{2}}$ for  $\delta <1/2$, then
$T_{\varphi}$ is a Toeplitz operator (for details see \cite{F}).
\end{rem}

Related with the definition of an operator $A$ by means of its
\emph{anti Wick symbol}, we can give an other representation of it
by means of the \emph{Wick symbol}: $\sigma^{W}(A)$, given by

\begin{eqnarray}\label{wickop}
Af(z)&=&\int e^{(z- v)\overline{v}} \, f(v) \,
\sigma^{W}(A)(\overline{v}, z) \,  dv(v)
\end{eqnarray}
where
\begin{equation}
\sigma^{W}(A)(\overline{v}, z)=\frac{\langle A K_{v}, K_{z}
\rangle }{ \langle K_{v}, K_{z} \rangle }.
\end{equation}

\begin{rem}
\ \\
 For all $f \in F^{2}(\CC^{n})$, we have that
$$f(v)=\langle f, K_{v}\rangle .$$
  Is in this sense that the coherent states are reproducing kernels of $F^{2}(\CC^{n})$.
   And then we can write
$$\sigma^{W}(A)(\overline{v}, z)=\frac{\langle A K_{v}, K_{z}
\rangle }{ \langle K_{v}, K_{z} \rangle }= \frac{ A K_{v}(z)}{K_{v}(z)}.$$

Now, taking into account that the norm of an operator can be bounded by its anti Wick and
Wick symbols:
$$
\|\sigma^{W}(A) \|_{\n} \leq \|A\| \leq \|\sigma^{AW}(A)\|_{\n},$$

we have that
 \begin{equation}\label{main}
 \frac{ |A K_{v}(z)|}{|K_{v}(z)|}\leq \|A\|.
\end{equation}

This relation between coherent states and Wick symbols will be needed in the proof of the main result.
\end{rem}

The correspondence between Wick or anti Wick symbols and
operators can be interpreted like a quantization process. Is in
this context that it is necessary to have a bilinear operation $*$
between the symbols, such that if $\tau= \varphi
* \psi$ then $T_{\tau}=T_{\varphi}T_{\psi}$.
This bilinear operation is completely defined for Wick symbols
obtaining a closed symbolic calculus. But the problem for anti
Wick symbols -as we have said in the Introduction- is still open.

\section {Some known results about composition of Toeplitz Operators.}

Among all the known positive results in  the composition of
Toeplitz operators on $F^2(\mathbb{C}^n)$, we have the following :
If $\varphi$ and $\psi$ are anti Wick symbols in one of these
classes
\begin{itemize}
    \item Analytic functions (\cite{Hal}) or
    \item Functions in the smooth Bochner algebra $B_{a}(\CC^{n})$
    (\cite{C}) or
    \item Polynomial functions  (\cite{C})
\end{itemize}

then $T_{\varphi}T_{\psi}=T_{\tau}=T_{\varphi  \diamond \psi}$,
where
\begin{equation}\label{diamond}
\varphi  \diamond \psi=
\sum_{k}\frac{(-1)^{k}}{k!}(\partial^{k}\varphi)(\overline{\partial}^{k}\psi).
\end{equation}

But unfortunately it is not possible to have this product in other
cases. In \cite{C} we can find an example that shows the
limitation on the ability to compose Toeplitz operators taking the
one dimensional Segal-Bargmann space $F^{2}(\CC)$ and
$\varphi=e^{2 \lambda|z|^{2}}$ with $\lambda=\frac{1+2i}{5}$.
Although $T_{\varphi}T_{\varphi}$ exists and, moreover, is a
bounded operator, there is not a symbol $\tau$ such that
$T_{\varphi}T_{\varphi}=T_{\tau}$ on $F^{2}(\CC)$.

In the next section, we will give a rather general result
encompassing this example, and a positive result about the
composition of Toeplitz operators with radial symbols in
$F^{2}(\CC)$.
\ \\

The theory about Toeplitz operators becomes interesting when the
symbols are radials, because these operators are unitary
equivalent to multiplication operators. More accurately, given a
radial symbol $a(|z|)=a(r)$, the Toeplitz operator $T_{a}$ is
unitary equivalent to the multiplication operator $\gamma_{a}I$
acting on $l_{2}^{+}$ where
\begin{equation}\label{gamma}
\gamma_{a}(n)=\frac{1}{n!}\int_{\RR^{+}}a(\sqrt{r})\, r^{n} \,
e^{-r}dr=\frac{2}{n!}\int_{\RR^{+}}a(r) \, r^{2n+1} \,
e^{-r^{2}}dr,
\end{equation}

i.e, there exist unitary operators $R: L^{2}(\CC, d\mu)\rightarrow
l_{2}^{+}$ and $R^{*}:l_{2}^{+}\rightarrow F^{2}(\CC)$ such that
\begin{equation}\label{RgammaR}
RT_{a}R^{*}\{c_{n}\}=\{\gamma_{a}(n).c_{n}\}
\end{equation}
for all sequence $\{c_{n}\} \in l_{2}^{+}$, where
\begin{equation}\label{operadorRstar}
R\varphi= \bigg\{\frac{1}{\sqrt{n!}}\int_{\CC}\varphi(z) \,
\overline{z}^{n} \,  d\mu (z)\bigg\}_{n \in \ZZ_{+}}
\end{equation}
and
$$R^{*}\{c_{n}\}= \sum_{n \in
\ZZ_{+}}\frac{c_{n}}{\sqrt{n!}} z^{n}.$$

Note that $RR^{*}=I:l_{2}^{+}\rightarrow l_{2}^{+}$ and
$R^{*}R=P:L^{2}(\CC, d\mu)\rightarrow {F}^{2}(\CC)$, where $P$ is
the projection defined in (\ref{projector}).

\begin{rem}
\ \\
In order to guarantee the convergence of the integral
(\ref{gamma}) it is considered, for example, $a(r)$ belonging to
$L_{1}^{\n}(\RR_{+},e^{-r^{2}})$. This space consists of all
measurable functions $a(r)$ on $\RR_{+}$ such that
    $$\int_{\RR^{+}} |a(r)| \, e^{-r^{2}} \, r^{n}dr < \n \ \ \mbox{for all $n \in \ZZ_{+}$.}$$

\end{rem}

Now, given a sequence $\gamma=\{\gamma(n)\} \in l_{\n}$
there exists a symbol $a(r)\in
L_{1}^{\n}(\RR_{+},e^{-r^{2}})$ such that the operator $T_{a}$ is
unitarily equivalent to the multiplication operator  by the
sequence $\gamma$, i.e $\gamma(n)=\gamma_{a}(n)$.

\begin{rem}\label{prop gamma}
\ \\
For $a(r)\in L_{1}^{\n}(\RR_{+},e^{-r^{2}})$ we have that
\begin{itemize}
    \item $T_{a}$ is a bounded operator if and only if
    $\{\gamma_{a}(n)\}$ is a bounded sequence.
    \item If $T_{a}$ is bounded, then its spectrum is given by
    $$sp (T_{a})=\overline{ \{ \gamma_{a}(n), n \in \ZZ_{+}\}}.$$
\end{itemize}

\end{rem}

The details of these results for radial symbols can be founded in
\cite{G-V}.
\ \\

\section{Main result.}

Let us define $L_{2}^{\n}(\RR_{+},e^{-r^{2}})$, the linear space
of all measurable functions $f(r)$ on $\RR_{+}$ such that
    $$q_{f}(n):= \int_{\RR^{+}} |f(r)|^{2} e^{-r^{2}}r^{n+1}dr < \n$$
for all $n \in \ZZ_{+}$. For this functions we define
$$A_{f}(x)=\sum_{n=0}^{\n}\frac{x^{n}q_{f}(n)}{n!}, \ \ x\in \RR.$$

\begin{thm}\label{main-th}
\ \\
Let  $T_{\varphi}$ and $T_{\psi}$ be two Toeplitz operators with
radial symbols $\varphi(r)$, $\psi(r) \in
L_{1}^{\n}(\RR_{+},e^{-r^{2}}))$ such that
\begin{enumerate}
\item $T_{\psi}$ is a bounded Toeplitz operator,
    \item  $\{ \gamma_{\varphi}(n). \gamma_{\psi}(n)\} \in
    l_{\n}$,
    \item $\varphi $ belongs to $
L_{2}^{\n}(\RR_{+},e^{-r^{2}})$ and $A_{\varphi}(x)$ converges
for all $x \in \RR_{+}$,
\end{enumerate}
then there exists $\tau \in L_{1}^{\n}(\RR_{+},e^{-r^{2}})$ such
that $T_{\tau}=T_{\varphi}T_{\psi}$ is a Toeplitz operator.
\end{thm}
\ \\
\begin{proof}
As we have mentioned in the previous section, if
$\{\gamma_{\varphi}(n). \gamma_{\psi}(n)\} \in
    l_{\n}$ there exists
$\tau \in L_{1}^{\n}(\RR_{+},e^{-r^{2}})$ such that
$\{\gamma_{\tau}(n)\}=\{\gamma_{\varphi}(n). \gamma_{\psi}(n)\}$.

Now, we have that:
\begin{enumerate}
    \item $T_{\tau}=T_{\varphi}T_{\psi}$ in the domain of $T_{\varphi}T_{\psi}$
    \item $T_{\tau}=T_{\varphi}T_{\psi}$ is indeed a Toeplitz
    operator.
\end{enumerate}

(1) Let us take $f(z) \in \mbox{\textit{Dom}}(T_{\varphi}T_{\psi})
\subset F^{2}(\CC)$. Then
$$f(z)=\sum_{n=0}^{\n}a_{n}z^{n}=R^{*}\{c_{n}\} ,$$
where $\{c_{n}\}=\{a_{n}\sqrt{n!}\}\in
l_{2}^{+}$ and $R^{*}$ is the operator given by (\ref{operadorRstar}).\\
From the definition of Toeplitz operator through the projection
$P=R^{*}R$ given in (\ref{toepl}) and by (\ref{RgammaR}) we have
that
\begin{eqnarray*}
T_{\varphi}T_{\psi}f &=& T_{\varphi}P(\psi f) \\
%   &=& T_{\varphi}P[\ \psi R^{*}\{c_{n}\}] \\
%   &=&  T_{\varphi}R^{*}R[\ \psi R^{*}\{c_{n}\}] \\
   &=&  T_{\varphi}R^{*}(R\psi R^{*})\{c_{n}\}  \\
   &=&  T_{\varphi}R^{*}\{\gamma_{\psi}(n). c_{n}\}\\
   &=&  R^{*}R\varphi R^{*}\{\gamma_{\psi}(n). c_{n}\}\\
     &=&  R^{*}\{\gamma_{\varphi}(n).\gamma_{\psi}(n)c_{n}\}=: R^{*}\{\gamma_{\tau}(n).c_{n}\}\\
     &=&R^{*}(RT_{\tau}R^{*})\{c_{n}\}:=T_{\tau}f
\end{eqnarray*}

(2) By the definition of Toeplitz operators, we must see that all the coherent states belong to the domain of
$T_{\tau}$, i.e $K_{a}(z)=e^{z\overline{a}} \in
\mbox{\textit{Dom}}(T_{\tau})$.

 Obviously $K_{a}(z) \in
\mbox{\textit{Dom}}(T_{\psi})$ since $T_{\psi}$ is a Toeplitz
operator. Now, it remains to prove that $\varphi T_{\psi} K_{a}(z) \in
L^{2}(\CC, d\mu)$. By the estimate (\ref{main}),

\begin{eqnarray*}
\int_{\CC} |\varphi(z) T_{\psi} K_{a}(z)|^{2} \, e^{-|z|^{2}} \,
dv(z)&\leq&
 \int_{\CC} |\varphi(z)|^{2} \, \|T_{\psi}\|^{2} \, |K_{a}(z)|^{2} \, e^{-|z|^{2}} \, dv(z) \\
\ \\
&\leq& \|{T_{\psi}}\|^{2} \, \int_{\CC} |\varphi|^{2} \,
\sum_{n=0}^{\n}\frac{(2|a||z|)^{n}}{n!} \, e^{-|z|^{2}} \, dv(z).
\end{eqnarray*}
Using polar coordinates, we have that
\begin{eqnarray*}
\|{T_{\psi}}\|^{2} \, \sum_{n=0}^{\n}\frac{|2a|^{n}}{n!}\int_{\CC}
|\varphi|^{2} \, |z|^{n}e^{-|z|^{2}}dv(z) &=& 2 \pi
\|{T_{\psi}}\|^{2}
 \sum_{n=0}^{\n}\frac{|2a|^{n}}{n!}\int_{\RR_{+}}
|\varphi|^{2} r^{n+1}e^{-r^{2}}dr\\
\ \\
&=& 2 \pi \|{T_{\psi}}\|^{2}  \sum_{n=0}^{\n}\frac{|2a|^{n}}{n!}q_{\varphi }(n)\\
\ \\
&=&\ 2 \pi \|{T_{\psi}}\|^{2} A_{\varphi }(2|a|)<\n\, ,
\end{eqnarray*}
since $\|{T_{\psi}}\|^{2} $ is bounded and  $A_{\varphi }(2|a|)$
exists for all $a \in \CC$.
%$$\int_{\CC} |\varphi T_{\psi} K_{a}(z)|^{2}e^{-|z|^{2}}dv(z)< \n .$$

Then $K_{a} \in \mbox{\textit{Dom}}(T_{\varphi}T_{\psi}) \subset
\mbox{\textit{Dom}}(T_{\tau})$, and therefore $T_{\tau}$ is indeed
a Toeplitz operator.

\end{proof}
\begin{flushright}
$\square$
\end{flushright}

\begin{rem}
\ \\
\textbf{Smooth Bochner Algebra.} For radial symbols $\varphi$ and
$\psi$ in the Bochner algebra $B_{a}(\CC)$, we know that
$T_{\varphi}T_{\psi}=T_{\varphi \diamond \psi}$
 with $\varphi \diamond \psi$ given by (\ref{diamond}) and
  $T_{\varphi \diamond \psi}$ is a Toeplitz operator. On the other hand, it is possible to
  apply the Theorem \ref{main-th}, and it can be proved that the symbol
  $\tau$ so  obtained satisfies $\tau = \varphi \diamond
    \psi$. We give the proof of these facts in Section 5.

\end{rem}

\begin{rem}
\ \\
As we have mentioned in Section I, there exist results about
the composition of Toeplitz operators for analytic symbols and
polynomial symbols. The first case is not contemplated by this
theorem, because radial symbols cannot be analytic (except for the constants where the theorem is trivial).\\
For polynomial radial symbols the situation is quite different:
When we calculate its sequence $\{\gamma(n)\}$ we have that it is
unbounded, then we cannot apply the Theorem \ref{main-th}.
Nevertheless the composition can be made using formula
(\ref{diamond}).
\end{rem}

\subsection{{\bf About the composition obstruction.}}

\begin{thm}\label{notcompose}
\ \\
 Let $A$ be a bounded operator with domain
in the Fock space $F^{2}(\CC)$, such that its Wick symbol is
$\sigma^{W}(A)=\, e^{-\theta|z|^{2}}$, with $ \theta \in \CC$.
Then
\begin{enumerate}
    \item  If $|\theta|^{2}= 2 Re(\theta)$ and $Re(\theta)>1$, there
is no bounded Toeplitz operator $T_{\tau}$ such that $A=T_{\tau}$.
    \item  If $|\theta|^{2} > 2 Re(\theta)$, there
is no bounded Toeplitz operator $T_{\tau}$ with $\tau \in
L_{1}^{\n}(\RR_{+},e^{-r^{2}})$ such that $A=T_{\tau}$.
\end{enumerate}
\end{thm}

\begin{proof}

(1) Let us take the operator $A$ such that
    $\sigma^{W}(A)=e^{- \theta|z|^{2}}$. Then by  (\ref{wickop})
\begin{eqnarray*}
Af(z)&=& \int e^{(z-
\xi)\overline{\xi}} \, f(v) \, \sigma^{W}(A)(z,\overline{\xi}) \, dv(\xi)\\
&=& \int e^{(z-\xi)\overline{\xi}} \, e^{- \theta(z
\overline{\xi})^{2}} \, f(v) \, dv(\xi).
\end{eqnarray*}

 If we define $e_{n}(z):=\frac{z^{n}}{\sqrt{n!}}$ , $\{e_{n}(z), n \in \NN \}$ is a basis of $F^{2}(\CC)$ and
$$Ae_{n}(z)= (1- \theta)^{n+1}e_{n}(z).$$

On the other hand, if we consider the operator
$\mathcal{M}_{a}f(z)=a f(az)$, we note that
$$Ae_{n}(z)=(1- \theta)^{n+1}e_{n}(z)=\mathcal{M}_{a}e_{n}(z)$$
where $a=1- \theta$. Then the operators $A$ and $\mathcal{M}_{a}$
agree in the basis and
 consequently in all the Fock space.\\

By the hypothesis on $ \theta$, we have that $|a|=|1-
\theta|=\sqrt{1-2 Re( \theta)+ | \theta|^{2}}=1$ and $Re (a)=1- Re
(\theta )<0$. By the results in \cite{C}, there is no bounded
Toeplitz operator $T_{\tau}$ such that $A=T_{\tau}$ in all the
Fock space.
\ \\

(2) Suppose that there exists $\tau \in L_{1}^{\n}(\RR_{+},e^{-r^{2}})$ such that $T_{\tau}$ is a bounded Toeplitz operator
    and  $T_{\tau} =A$.

   The Wick symbol of $T_{\tau}$ can be calculated by
\begin{equation}\label{wickgamma}
\sigma^{W}(T_{\tau})=e^{-|z|^{2}}\sum_{n=0}^{\n}\frac{|z|^{2n}}{n!}\gamma_{\tau}(n)
\end{equation}
where $\gamma_{\tau}(n)$ is defined by (\ref{gamma}) (see
\cite{G-V} for more details).

On the other hand
\begin{equation}
e^{-
\theta|z|^{2}}=e^{-|z|^{2}}\sum_{n=0}^{\n}\frac{|z|^{2n}}{n!}(1-
\theta)^{n}.
\end{equation}
Then, taking into account that
$$\sigma^{W}(T_{\tau})= e^{- \theta|z|^{2}}$$

we have that
$$|\gamma_{\tau}(n)|= |(1- \theta)^{n}|=({1-2 Re( \theta)+ |\theta|^{2}})^{\frac n2}\geq
(1+\epsilon)^{\frac n2}.$$ Then $\gamma_{\tau}(n)$ is an unbounded
sequence. This contradicts that the operator can be bounded as it was
stated in Remark \ref{prop gamma}.
\end{proof}
\begin{flushright}
$\square$
\end{flushright}

\begin{rem}
\ \\
This theorem  includes the example mentioned in Section 3 given in \cite{C}:\\
If $\varphi=e^{2(\frac{1+2i}{5})|z|^{2}}$,
$$\gamma_{\varphi}(n)= \frac{1}{n!}\int_{\RR^{+}}e^{2(\frac{1+2i}{5}) r}r^{n}e^{-r}dr
= \left( \frac{3}{5}-\frac{4}{5}i \right) ^{-(n+1)}\, ,$$

and then $\{\gamma_{\varphi}(n)\} \in l_{\n}$ ($|\gamma_{\varphi}(n)|=1$). The Wick symbol of $T_{\varphi}T_{\varphi}$ can be calculated by (\ref{wickgamma}),
$$ \sigma^{W}(T_{\varphi}T_{\varphi})(z,\overline{z})=
e^{-|z|^{2}}\sum_{n=0}^{\n}\frac{|z|^{2n}}{n!}\gamma_{\varphi}(n)\gamma_{\varphi}(n)= C e^{-K|z|^{2}}\, $$

with $C=\left(\frac{3}{5}-\frac{4}{5}i \right)^{-2}$ and $K=\left(\frac{32}{25}+\frac{24}{25}i\right)$. Then $|K|^{2}=2Re(K)=\frac{64}{25}$ and by (1) in Theorem \ref{notcompose}, there is not a bounded Toeplitz operator with this symbol, and as $T_{\varphi}T_{\varphi}$ is a bounded operator, it cannot be a Toeplitz operator on $F^{2}(\CC)$.
\end{rem}

\section{Smooth Bochner Algebra.}

Consider $\varphi(|z|)=\varphi(r)$, $\psi(|z|)=\psi(r)$ radial
functions in the Bochner algebra $B_{a}(\CC)$. The functions in
$B_{a}(\CC)$  are bounded, uniformly continuous, with bounded
derivatives of all orders. Therefore, by the hypothesis on
$\varphi$ and $\psi$ we can apply the Theorem \ref{main-th} and
then, there exists $\tau$ such that $T_{\tau}=T_{\varphi}T_{\psi}$
is a Toeplitz operator in the Fock space.

As we have mentioned in  Section 3, we know that
$T_{\varphi}T_{\psi}=T_{\varphi \diamond \psi}$
 with $\varphi \diamond \psi$ given by (\ref{diamond})
 and $T_{\varphi \diamond \psi}$ is a Toeplitz operator.

Then we have on $F^{2}(\CC)$ that $T_{\tau}=T_{\varphi \diamond
\psi}$. The question is if $\tau = \varphi \diamond \psi$ since
for Toeplitz operators, is not always true that if $T_{\eta}=0$
then $\eta=0$.

However, we can prove that if the anti Wick symbol $\eta$ of
$T_{\eta}$ belongs to $E^{\varepsilon}(\RR_{+},e^{-r^{2}})$ (a
subclass of $L_{1}^{\n}(\RR_{+},e^{-r^{2}})$), then the operator
$T_{\eta}=0$ if and only if $\eta=0$ almost everywhere (for
details see \cite{G-V}).

The space $E^{\varepsilon}(\RR_{+},e^{-r^{2}})$ is the subclass of
$L_{1}^{\n}(\RR_{+},e^{-r^{2}})$ that consists of all functions
$\eta(r)$ satisfying at $+\n$ the following estimate: fixed
$\varepsilon >0$,
$$|\eta(r)|e^{-r^{2}+\varepsilon r}\leq C.$$

\begin{prop}
\ \\
Given $\varphi$ and $\psi$ two radial functions in the Bochner
algebra $B_{a}(\CC)$. Then
\begin{enumerate}
    \item $\varphi$ and $\psi$ satisfies the hypothesis of Theorem
\ref{main-th} and
    \item the function given by Theorem \ref{main-th} satisfies $\tau = \varphi \diamond
    \psi$.
\end{enumerate}
\end{prop}

\begin{proof}
We must prove that there exists $\varepsilon
> 0$ such that $ \tau - \varphi \diamond \psi \in
E^{\varepsilon}(\RR_{+},e^{-r^{2}})$.

 $\underline{ \tau \in E^{\varepsilon}(\RR_{+},e^{-r^{2}})}:$
\begin{eqnarray*}
\int_{1}^{\n}|\tau(r)|e^{-r^{2}+\varepsilon r} dr
 &\leq& \left( \int_{1}^{\n}|\tau(r)|^{2} r e^{-r^{2}}dr \right)^{1/2} \,
 \left( \int_{1}^{\n} \frac{e^{-r^{2}}e^{2 \varepsilon
 r}}{r}dr\right)^{1/2}\\
 \ \\
  &\leq& \left( \int_{1}^{\n}|\tau(r)|^{2} r e^{-r^{2}}dr \right)^{1/2} \,
 \left( \int_{1}^{\n}e^{-r^{2}}e^{2 \varepsilon r}dr\right)^{1/2}.
\end{eqnarray*}

It can be easily seen  that the second integral in the product converges for all
$\varepsilon $.

Also, note that $1 \in F^{2}(\CC)$ and $T_{\tau}$ is bounded, then
$1 \in \mbox{\textit{Dom}}(T_{\tau})$, i.e.,
$$\int_{\CC}|\tau . 1|^{2}e^{-|z|^{2}} dv(z)=(2 \pi)\int_{0}^{\n}|\tau(r)|^{2}r e^{-r^{2}}dr < \n.$$
Then the first integral converge.\\

Therefore $\int_{1}^{\n}|\tau(r)|e^{-r^{2}+\varepsilon r}dr $
converge and then  $|\tau(r)|e^{-r^{2}+\varepsilon r}$ is bounded
in $+ \n$.

Taking into account that $\varphi \diamond \psi$ is bounded, then
$\eta(r):=\tau(r)-\varphi \diamond \psi(r) \in
E^{\varepsilon}(\RR_{+},e^{-r^{2}})$.
\begin{flushright}
$\square$
\end{flushright}
\end{proof}
\section{Appendix.}

A {\bf Fock space} $F$ is a Hilbert space in which there exist two
operators $\widehat{a}$ and $\widehat{a}^{*}$ called annihilation
and creation operators, satisfying the canonical commutation
rules, i.e $[\widehat{a},\widehat{a}^{*}]=I$, where $I$ is the
identity operator. There exists also a vector $\Phi_{0}$ (called
\emph{ vacuum vector}) annihilated by the annihilations operators,
such that the system
$$\left\{\frac{((\widehat{a}*)^{\alpha}\Phi_{0})}{\sqrt{\alpha
!}}\right\}$$ is complete and orthonormal in $F$.

Classical examples of Fock spaces are:

\begin{itemize}
    \item $L^{2}(\RR^{n})$: The Hilbert space of square
integrable functions on $\textsl{R}^{n}$, with the Lebesgue
 measure on $\textsl{R}^{n}$.

    \item $F^{2}(\CC^{n}, d\mu)$: the subspace of $L^{2}(\CC^{n}, d\mu)$ of all
the analytic functions in $\CC^{n}$. Here
$d\mu(z)=\pi^{-n}e^{-|z|^{2}}dv(z)$ is the Gaussian measure, where
$dv(z)$ is the Lebesgue measure on $\textsl{C}^{n}$.
\end{itemize}

In the next table, we show the annihilation and creation
operators, the vacuum vector and the basis of each space.

\begin{center}
\begin{tabular}{|c|c|c|}
  \hline
  % after \\: \hline or \cline{col1-col2} \cline{col3-col4} ...
   & $L^{2}(\RR^{n})$ & $F^{2}(\CC^{n}, d\mu)$ \\
  \hline
 % &&\\
  Creation operator $\widehat{a}^{*}$ & $\frac{\widehat{q}-i\widehat{p}}{\sqrt{2}}$ & $\widehat{a}^{*}f(z)=z.f(z)$ \\
 % &&\\
 Annihilation operator $\widehat{a}$ & $\frac{\widehat{q}+i\widehat{p}}{\sqrt{2}}$ & $\widehat{a}f(z)=\frac{\partial f(z)}{\partial z}$ \\
% &&\\
 Vacuum vector & $\Phi_{0}(x)=\pi^{-n/4}e^{-x^{2}/2}$ & $\Phi_{0}(z)=1$ \\
%&&\\
 Basis  & Hermite functions &
$\left\{z^{\alpha}/ \sqrt{\alpha
!}\right\}$ \\
%&&\\
  \hline
\end{tabular}
\end{center}

There is a unitary isomorphism that apply a basis of $L^{2}(\RR)$
onto a basis of $F^{2}(\CC^{n}, d\mu)$. This map is called the
\emph{Bargmann transform} and is given by
 $$L\varphi (\overline{z})=(\pi h)^{-(n/4)}\int
e^{-\frac{1}{2}(s^{2}-2\sqrt{2}s.\overline{z}+\overline{z}^{2})}\varphi(s)ds\,
.$$ This transform allows us to interpret an operator in
$F^{2}(\CC^{n}, d\mu)$ as a pseudodifferential operator in
$L^{2}(\RR^{n})$.

The coherent states (or Poisson vectors) are the functions
$K_{v}(z)=e^{z\overline{v}} \in F^{2}(\CC^{n}, d\mu)$, indeed they
are the eigenfunctions of the operator
$\widehat{a}=\frac{\partial}{\partial z}$ associated to the
eigenvalue $\overline{v}$. By the Bargmann transform, we can
obtain that $f(v)=(f, K_{v})$ for all $f \in F^{2}(\CC^{n},
d\mu)$. In this sense, $K_{v}$ is a reproductive Kernel and
consequently $ F^{2}(\CC^{n}, d\mu)$ is a Reproducing Kernel
Hilbert Space (RKHS).

The Wick and anti Wick operators appear naturally in the context
of quantization procedures, the correspondence between operators
and symbols is related with the association between a quantum and
a classical observable. There are different ways to do it
consistent with the probabilistic interpretation by mean of
pseudodifferential operators.

In $F^{2}(\CC^{n}, d\mu)$, taking into account that the
composition of the annihilation and construction operator is not
commutative, we consider the operator theory associated to the
product order
    $\widehat{a}\widehat{a}^{*}$ (\emph{Wick}),
    or $\widehat{a}^{*}\widehat{a}$ (\emph{Anti Wick}).

In $L^{2}(\RR^{n})$ the Weyl calculus gives the suitable frame in
quantization procedures.

 All these symbols are related by the
kernel of the heat $H_{t}$ in the following way, where
$\{H_{t}\}_{t>0}$ on $\RR^{2n}$ is defined by $H_{t}f= f*
\gamma_{t}$ with $\gamma_{t}(x, \xi)=(4t)^{-n}e^{-(x^{2}+
\xi^{2})/4t}$:
\begin{equation}\label{heat1}
 \sigma_{w}(A)=
H_{1/2}(\sigma^{AW}(A))\,\,\, , \sigma^{W}(A)=
H_{1/2}(\sigma_{w}(A))\,\,\, ,  \sigma^{W}(A)=
H_{1}(\sigma^{AW}(A)).
\end{equation}
(See for instance \cite{S-B} or \cite{F}.)

The Weyl and Wick symbols of a given pseudodifferential operator
are always defined. By (\ref{heat1}) we can see that \textbf{not
all} these operators have an anti Wick symbol, since in order to
find an anti Wick symbol one has to solve -for example- the
inverse heat equation with initial condition given by Weyl symbol
for the time $t=1/2$.

Nevertheless, the symbol anti Wick have two
advantages that the other symbols do not verify at the same time:
\begin{enumerate}
    \item \textbf{Positivity:} If $\varphi \geq 0$ then $T_{\varphi}f \geq
    0$ for all $f$ in the domain of $T_{\varphi}$.
    \item \textbf{Selfadjointness:} If $\varphi$ is a function to real values then $T_{\varphi}$ is selfadjoint.
\end{enumerate}

\smallskip

\end{document}